\documentclass[11pt,fleqn]{article}

\usepackage[cp1251]{inputenc}

\usepackage{amsmath,amssymb,amsthm,enumerate,epsfig,graphics,cite}
\usepackage[ukrainian,english]{babel}

\newcommand{\todo}[1][\null]{\ensuremath{\clubsuit}}

\newcommand{\noprint}[1]{}

\begin{document}
\begin{center}
\Large\bf
New generalizations of Zeta-Function and Tricomi function
\end{center}
\begin{center}
\it
N. Virchenko, A. Ponomarenko
\end{center}

{\it

We present new generalizations of the zeta function, Tricomi functions. Their basic properties are studied.
These new generalizations are carried out with the help of $(\tau,\beta)$-the generalized confluent hypergeometric function.}


In connection with the wide use of special functions in probability theory and mathematical statistics, continuum mechanics, quantum mechanics, quantum optics, astrophysics, diffraction theory, aerodynamics, coding theory, biomedicine and other applications of special functions, the interest to them increased significantly \cite{Andrews1985,Andrews1999,Edwards1954,Ewell1990,Titchmarsh1951}. In particular, this can be said about the
 generalized hypergeometric Wright functions has increased \cite{Wright1940}.
Wright introduced a generalization of the hypergeometric function in the form of a series \cite{Wright1940}
\begin{gather}\label{eq1}
\sum_{n=0}^\infty\frac{\Gamma(\alpha_1+\tau_1 n)\cdots \Gamma(\alpha_p+\tau_p n)z^n}{\Gamma(\rho_1+\beta_1 n)\cdots \Gamma(\rho_q+\beta_q n)n!},\quad z\in {\mathbb C},
\end{gather}
in which the parameters $\tau_i, \beta_j \in {\mathbb R}_+$ $(i= \overline{1,p}$, $j= \overline{1,q}) $ satisfy the condition
\begin{gather}\label{eq2}
 1+\sum_{k=0}^{g} \beta_k - \sum_{m=1}^{p} \tau_m >0.
\end{gather}

Under the condition \eqref{eq2}, the series \eqref{eq1} converges for any $ z $. Under the condition
\begin{gather}\label{eq3}
 1+\sum_{k=0}^{g} \beta_k - \sum_{m=1}^{p} \tau_m =0,
\end{gather}
the series converges for any $ |z|< \dfrac{1}{r} $.

In this paper we introduce new generalizations of the zeta function, the Tricomi functions;
their main properties are studied. This opens the way to a deeper, better application of
these functions both in the theory of special functions, and in many applied sciences.
These new generalizations are carried out with the help of $(\tau, \beta)$-confluent hypergeometric  function \cite{Virchenko2006}:
\begin{gather}\label{eq4}
{} ^{r} _1 \Phi_1^{\tau, \beta} (a;c;z)\equiv {}^{r} _1 \Phi_1^{\tau, \beta} (z)= \frac{\Gamma(c)}{\Gamma(a)\Gamma(c-a)} \int _{0}^{1} t^{a-1} (1-t) ^{c-a-1} {} _1  \Psi _1 \left[\begin{array}{cc} (c;\tau); & \\
 &  zt^\tau\\
 (c;\beta); & \end{array} \right] dt,
\end{gather}
where $ \operatorname{Re} c > \operatorname{Re} a > 0$, $\{\tau,\beta\}\subset{\mathbb R}_+$; $\tau >0$; $\beta >0$, $\tau - \beta < 1$, $\Gamma(a)$ is the classical gamma function \cite{Beitmen1965}, $ {}_1 \Psi_1 $ is the Fox--Wright function \cite{Wright1940}:
\begin{gather}
{} _1  \Psi _1 \left[\begin{array}{cc} (a;\alpha); & \\
 &  z\\
 (b;\beta); & \end{array} \right] = \sum_{k=0}^\infty \frac{\Gamma(a+\alpha k)}{\Gamma(b+\beta k)}\frac{z^k}{k!}.  \nonumber
\end{gather}

Earlier we considered \cite{Virchenko2016} generalizations of the gamma, beta, psi, Laguerre functions, and
Volterra functions.

\textbf{1. The generalized zeta-function.} As it is known \cite{Beitmen1965}, the Riemann zeta-function of the form
\begin{gather}\label{eq5}
\zeta (\alpha)= \sum_{n=1}^\infty\frac{1}{n^\alpha},
\end{gather}
where $\alpha=\sigma+i\tau$, $\sigma>1$, was considered by Euler in 1737 for real $\alpha $. He also proved the identity
\begin{gather}\label{eq6}
\zeta (\alpha)= \prod_{p} \left(1- \frac{1}{p^\alpha}\right)^{-1}, \quad \alpha=\sigma + i \tau, \quad\sigma>1.
\end{gather}
We introduce a generalization of the Riemann zeta-function in the form
\begin{gather}\label{eq7}
\zeta ^r (\alpha)= \frac {1}{\Gamma(\alpha)} \int _{0}^\infty \frac{t^{\alpha-1}e^{-t}}{1-e^{-t}} {}^{r} _1 \Phi_1^{\tau,\beta} \left(a;c;-\frac{r}{t}\right) {\rm d}t,
\end{gather}
where $\operatorname{Re} r > 0$; $r=0$, $\sigma >1 $.

Let us study the basic properties of $\zeta ^r (\alpha)$.

\textbf{Theorem 1.}   {\it We have the equality}
\begin{gather}\label{eq8}
\zeta ^r (\alpha)= \frac {1}{\Gamma(\alpha)} \sum_{n=1}^\infty \Gamma_{nr} (\alpha) n^{-\alpha}, \quad \operatorname{Re} r > 0;\quad r=0, \quad\sigma >1 ,
\end{gather}
{\it where} $\Gamma_{nr}$ {\it is the generalized gamma function \cite{Virchenko2016}:}
\begin{gather}
\Gamma_{\varsigma} (\alpha)=\int _{0}^\infty t^{\alpha-1}e^{-t} \, {}^{r} _1 \Phi_1^{\tau, \beta} \left(a;c;-\frac{\varsigma}{t}\right){\rm d}t, \quad \operatorname{Re} c > \operatorname{Re} a > 0,\operatorname{Re} \alpha > 0, \varsigma >0, \{\tau,\beta\} \in {\mathbb R}_+; \nonumber \\  \tau >0; \beta >0, \tau - \beta < 1. \nonumber
\end{gather}

\textbf{Proof.} Expanding the function  $\left(1-e^{-t}\right)^{-1}$ in a series in powers of $ e^{-t} $, we get
\begin{gather}\label{eq9}
\zeta ^r (\alpha)= \frac {1}{\Gamma(\alpha)} \int_{0}^\infty \sum_{n=1}^\infty t^{\alpha-1} e^{-nt} \,{}^{r} _1 \Phi_1^{\tau,\beta} \left(a;c;-\frac{r}{t}\right) {\rm d}t.
\end{gather}

We change the order of summation and integration in \eqref{eq9} (this is possible, since the series
converges to an integrable function). Further, after the transformations, taking into account the definition of the generalized gamma function, we obtain  \eqref{eq8}.

\textbf{Notes.}

1. Using the property of the generalized gamma function \cite{Virchenko2016}, we obtain the equality:
\begin{gather}\label{eq10}
\sum_{n=1}^\infty\ \Gamma_{nr}(\alpha)=r^{\alpha}\Gamma(-\alpha)\zeta^{r}(-\alpha),\quad \operatorname{Re} r > 0, \quad \alpha\neq 0,1,2,....
\end{gather}

2. After using the equality \cite{Beitmen1965}
\begin{gather}\label{eq11}
\Gamma_{r}(\alpha)=2r^{\frac{\alpha}{2}}K_{\alpha}\left(2\sqrt{r}\right),  \quad \operatorname{Re} r > 0, \quad \left|\operatorname{arg}\sqrt{r}\right|<\pi,
\end{gather}
we get the connection between the generalized zeta function and the MacDonald function:
\begin{gather}\label{eq12}
\zeta^{r}(\alpha)=\frac{2r^{\frac{\alpha}{2}}}{\Gamma(\alpha)} \sum_{n=1}^\infty n^{-\frac{\alpha}{2}}K_\alpha \left(2\sqrt{nr}\right),\quad\operatorname{Re} r > 0.
\end{gather}

3. After replacing $t$ by $2t$ with \eqref{eq7} and after using the equality
\begin{gather}\label{eq13}
\left(e^{2t}-1\right)^{-1}=e^{-t}\operatorname{csc}h(t)/2,
\end{gather}
we obtain an integral representation of the function $\zeta^{r}(\alpha)$:
\begin{gather}\label{eq14}
\zeta^{r}(\alpha)=\frac{2^{\alpha-1} }{\Gamma(\alpha)} \int_{0}^\infty t^{\alpha-1}e^{-t} \, {}^{r} _1 \Phi_1^{\tau,\beta} \left(a;c;-\frac{r}{2t}\right)\operatorname{csc}h(t) {\rm d}t.
\end{gather}
If $r = 0$ then we have a known result for the classical zeta-function:
\begin{gather}\label{eq15}
\zeta(\alpha)=\frac{2^{\alpha-1} }{\Gamma(\alpha)} \int_{0}^\infty t^{\alpha-1}\operatorname{csc}h(t) {\rm d}t, \quad \sigma>1.
\end{gather}

\textbf{Theorem 2.}  {\it For the generalized zeta-function we have the equality }
\begin{gather}\label{eq16}
\zeta^{r}(\alpha)-2^{1-\alpha}\zeta^{2r}(\alpha)=\frac{1}{\Gamma(\alpha)}\int_{0}^\infty \frac {t^{\alpha-1}e^{-t}} {1+e^{-t}}  {}^{r} _1 \Phi_1^{\tau,\beta} \left(a;c;-\frac{r}{t}\right) {\rm d}t,  \quad \operatorname{Re} r > 0; r=0,\sigma>1.
\end{gather}

\textbf{Proof.} From the formula for the generalized function $\Gamma$-function we obtain
\begin{gather}\label{eq17}
n^{-\alpha}\Gamma_{nr}(\alpha)=\int_{0}^\infty\ t^{\alpha-1}e^{-nt} \, {}^{r} _1 \Phi_1^{\tau,\beta} \left(a;c;-\frac{r}{t}\right) {\rm d}t.
\end{gather}
after multiplying by $\left(-1\right)^{n-1}$, and summing both sides of the previous equality by $n$, we get
\begin{gather}\label{eq18}
\sum_{n=1}^\infty\ \left(-1\right)^{n-1} \frac{\Gamma_{nr}(\alpha)}{n^{\alpha}}=\int_{0}^\infty t^{\alpha-1}\left(1+e^{-t}\right)^{-1}e^{-t} \, {}^{r} _1 \Phi_1^{\tau,\beta} \left(a;c;-\frac{r}{t}\right) {\rm d}t.
\end{gather}
now from the equality \eqref{eq18} and the formula
\begin{gather}\label{eq19}
\zeta^{r}(\alpha)-2^{1-\alpha}\zeta^{2r}(\alpha)=\frac{1}{\Gamma(\alpha)} \sum_{n=1}^\infty\ \left(-1\right)^{n-1} \frac{\Gamma_{nr}(\alpha)}{n^{\alpha}},\quad  \operatorname{Re} r > 0;\quad r=0,\quad \sigma>1,
\end{gather}
we get \eqref{eq16}.

\textbf{Corollary.} From the formula \eqref{eq16} with $r = 0$ we obtain the following {\it formula for the classical zeta-function:}
\begin{gather}\label{eq20}
\zeta(\alpha)=\frac{1}{\Gamma(\alpha)(1-2^{1-\alpha})}\int_{0}^\infty\ t^{\alpha-1}\left(1+e^t\right)^{-1}{\rm d}t, \quad \sigma>1.
\end{gather}

Of considerable interest is the relationship between the generalized zeta function and integral
transformations.

It is obvious that \eqref{eq7} can be written in the form of an integral Mellin transform:
\begin{gather}\label{eq21}
\zeta^{r}(\alpha)=\frac{1}{\Gamma(\alpha)}M \left \{\frac{e^{-t} \, {}^{r} _1 \Phi_1^{\tau,\beta} \left(a;c;-\frac{r}{t}\right)}{1-e^{-t}};\alpha \right\}, \quad \operatorname{Re} r > 0;\quad r=0,\quad \sigma>1.
\end{gather}
After performing the substitution $t = e^{-x}$, we obtain the integral Laplace transform:
\begin{gather}\label{eq22}
\zeta^{r}(\alpha)=\frac{1}{\Gamma(\alpha)}\int_{-\infty}^\infty\ \frac{ {}^{r} _1 \Phi_1^{\tau,\beta} \left(a;c;-re^x\right)}{exp(e^{-x})-1}e^{-\alpha x}{\rm d}x.
\end{gather}
Performing the substitution $\alpha = \sigma+ i\omega$ in \eqref{eq22} leads to an integral Fourier transform:
\begin{gather}\label{eq23}
\zeta^{r}(\alpha)=\frac{1}{\Gamma(\alpha)}\int_{-\infty}^\infty\ \frac{ {}^{r} _1 \Phi_1^{\tau,\beta} \left(a;c;-re^x\right)}{exp(e^{-x})-1}e^{-i\alpha\omega x}{\rm d}x.
\end{gather}

{\it The second generalization of the zeta function} can be represented in the form
\begin{gather}\label{eq24}
^* \zeta^{r}(\alpha)=\frac{1}{A(\alpha)}\int_{0}^\infty\ \frac { t^{\alpha-1} e^{-t}} {1+e^{-t}}  {}^{r} _1 \Phi_1^{\tau,\beta} \left(a;c;- \frac{r} {t}\right) {\rm d}t,
\end{gather}
where $\operatorname{Re} r > 0 $; $r=0$, $\sigma>1$, $A(\alpha)=\Gamma(\alpha)(1-2^{1-\alpha})$.

Representing the expression $\left(1+e^{-t}\right)^{-1}$  by a series, we get
\begin{gather}\label{eq25}
^* \zeta^{r}(\alpha)=\frac{1}{A(\alpha)} \sum_{n=1}^\infty\ \left(-1\right)^{n-1} \int_{0}^\infty\ t^{\alpha-1} e^{-nt}  \, {}^{r} _1 \Phi_1^{\tau,\beta} \left(a;c;- \frac{r} {t}\right) {\rm d}t,
\end{gather}
and after transformations we have
\begin{gather}\label{eq26}
 ^* \zeta^{r}(\alpha)=\frac{1}{A(\alpha)} \sum_{n=1}^\infty\ \frac {\left(-1\right)^{n-1}} {n^\alpha} \Gamma_{nr}(\alpha), \quad   \operatorname{Re} r > 0;\quad r=0,\quad \sigma>1.
\end{gather}

{\it A generalization of the Hurwitz zeta-function} \cite{Beitmen1965}
\begin{gather}\label{eq27}
\zeta(\alpha,q)= \sum_{n=0}^\infty\ \frac{1}{(n+q)^\alpha}=\frac{1}{\Gamma(\alpha)} \int_{0} ^\infty\ \frac{t^{\alpha-1} e^{-qt}}{1-e^{-t}} {\rm d}t,
\quad \sigma>1, \quad 0 < q \leq 1,
\end{gather}
is representable in the form
\begin{gather}\label{eq28}
\zeta^{r} (\alpha,q)= \sum_{n=0}^\infty\ \frac{ {}^{r} _1 \Phi_1^{\tau,\beta} \left(a;c;- \frac{r} {t}\right)}{\left(n+q\right)^{\alpha}}=\frac{1}{\Gamma(\alpha)} \int_{0} ^\infty\ \frac{t^{\alpha-1} e^{-qt}} {1-e^{-t}} {} ^{r} _1 \Phi_1^{\tau,\beta} \left(a;c;- \frac{r}{t}\right) {\rm d}t,
\end{gather}
\begin{gather}\label{eq29}
^* \zeta^{r} (\alpha,q)=\frac{1}{A(\alpha)}\int_{0}^\infty\ \frac{ t^{\alpha-1} e^{-qt}} {1-e^{-t}}  {}^{r} _1 \Phi_1^{\tau,\beta} \left(a;c;- \frac{r}{t}\right) {\rm d}t,
\end{gather}
where $0 < q \leq 1$, $\operatorname{Re} r > 0$; $r=0$, $\sigma>1$.

It is easy to see that
\begin{gather}\label{eq30}
\zeta (\alpha,q) = \zeta^{r} (\alpha,q)| _{r=0}, \quad \sigma>1 ,
\end{gather}
and
\begin{gather}\label{eq31}
^* \zeta^{r} (\alpha,q) = \frac {\zeta ^r (a,q)} {1-2^{1-\alpha}}, \quad \sigma>1, \quad 0< q\leq 1.
\end{gather}

\textbf{Theorem 3.}  {\it The following equality holds}
\begin{gather}\label{eq32}
^* \zeta^{r} (\alpha,q) = \frac {2^{1-\alpha} \zeta ^{2r} (a,\frac{q}{2}) - \zeta ^r (a,q)} {1-2^{1-\alpha}}, \quad \sigma>1, \quad 0< q\leq 1.
\end{gather}

\textbf{Proof.} We change $q$ to $ \frac{q}{2}$ in \eqref{eq28}, $t$ to $2t$, we get
\begin{gather}\label{eq33}
\zeta^{2r}\left(\alpha,\frac{q}{2}\right)=\frac{2^{\alpha}}{\Gamma(\alpha)}\int_{0}^\infty\ \frac{t^{\alpha-1} e^{-qt}} {1-e^{-2t}} {}^{r} _1 \Phi_1^{\tau,\beta} \left(a;c;- \frac{r}{2t}\right) {\rm d}t.
\end{gather}
We also have the equality
\begin{gather}\label{eq34}
2\left(1-e^{-2t}\right)=\left(1-e^{-t}\right)^{-1}+\left(1+e^{-t}\right)^{-1}.
\end{gather}

From the equations \eqref{eq33} and \eqref{eq34} we obtain
\begin{gather}
2^{1-\alpha}\zeta^{2r}\left(\alpha,\frac{q}{2}\right) - \zeta^{r}(\alpha,q) = \left(1-2^{1-\alpha}\right){}^*\zeta^{r}(\alpha,q), \nonumber
\end{gather}
and then we get \eqref{eq32}. Using the relation
\begin{gather}
\frac{1}{e^{t}+1} = \frac{1}{e^{t}-1} - \frac{2}{e^{2t}-1},  \quad  \int_{0}^\infty\ \frac{t^{\alpha-1}}{e^{2t}-1} {\rm d}t = 2^{-\alpha}\int_{0}^\infty\ \frac{t^{\alpha-1}}{e^{t}-1}{\rm d}t, \nonumber
\end{gather}
we have an analytic continuation
\begin{gather}\label{eq35}
\zeta^{r}(\alpha) = \frac{1}{\Gamma(\alpha)\left(1-2^{1-\alpha}\right)} \int_{0}^\infty\ \frac{t^{\alpha-1}}{e^{t}+1} {\rm d}t.
\end{gather}

\textbf{2. New generalization of Tricomi function.} Tricomi in 1927 introduced the function $\Psi (a, c; x)$ \cite{Tricomi1954}:
\begin{gather}\label{eq36}
\Psi (a, c; x) = \frac{1}{\Gamma(a)} \int_{0}^\infty\ t^{a-1}(1+t)^{c-a-1}e^{-xt}{\rm d}t, \quad \operatorname{Re} a > 0,
\end{gather}
which determines a solution of the equation \cite{Beitmen1965}
\begin{gather}\label{eq37}
x \frac{{\rm d}^2 y}{{\rm d} x^2}  +  (c-x)\frac{{\rm d} y}{{\rm d} x} -ay = 0
\end{gather}
in the field $\operatorname{Re} x > 0$. We can extend the domain of definition by rotating the path of integration.

We introduce a generalization of the Tricomi function \eqref{eq36} in the form
\begin{gather}\label{eq38}
^r U^{\tau,\beta} (a,c;\alpha;\gamma;\delta;x)=\frac{1}{\Gamma(a)} \int_{0}^\infty\ t^{a-1}(1+t)^{c-a-1}e^{-xt} \, {} ^{r} _1 \Phi_1^{\tau,\beta} \left(\alpha;\gamma;- \frac{r} {t^{\delta}}\right) {\rm d}t,
\end{gather}
where $\operatorname{Re} a >\operatorname{Re} c > 0$, $\{\tau,\beta\}\subset{\mathbb R}_+$; $\tau>0$; $\beta >0$, $\tau-\beta<1$, $\{a,c\}\in {\mathbb C}$, $\delta>0$, $r>0$, \\ $\operatorname{Re} \alpha >\operatorname{Re} \gamma > 0$, $ ^{r} _1\Phi_1^{\tau,\beta} (...)$ --- $(\tau, \beta)$-comfluent hypergeometric function, which is defined by the formula \eqref{eq4}, and can also be represented by the following series \cite{Virchenko2016}:
\begin{gather}\label{eq39}
{} ^{r} _1 \Phi_1^{\tau,\beta} \left(\alpha;\gamma;- \frac{r} {t^{\delta}}\right)=\frac{\Gamma(\gamma)}{\Gamma(\alpha)}\sum_{n=0}^\infty\ \frac{\Gamma(\alpha+\tau n)}{\Gamma(\gamma+\beta n)}\left (-\frac{r}{t^{\delta}} \right) ^n \frac{1}{n!}.
\end{gather}

Differentiation formulas for the function $ ^r U^{\tau,\beta} $ are obtained by direct differentiation:
\begin{gather}\label{eq40}
\frac{{\rm d} }{{\rm d} x}{} ^r U^{\tau,\beta} (a,c,\alpha;\gamma;\delta;x)=-\frac{1}{\Gamma(a)} \int_{0}^\infty\ t^{a}(1+t)^{c-a-1}e^{-xt} \, {} ^{r} _1 \Phi_1^{\tau,\beta} \left(\alpha;\gamma;- \frac{r} {t^{\delta}} \right) {\rm d}t,
\end{gather}
\begin{gather}\label{eq41}
\frac{{\rm d}^n }{{\rm d} x^n}{} ^r U^{\tau,\beta} (a,c;\alpha;\gamma;\delta;x)=\frac{(-1)^n}{\Gamma(a)} \int_{0}^\infty\ t^{a+n-1}(1+t)^{c-a-1}e^{-xt} \, {} ^{r} _1 \Phi_1^{\tau,\beta} \left(\alpha;\gamma;- \frac{r} {t^{\delta}}\right) {\rm d}t.
\end{gather}
The following theorem holds.

\textbf{Theorem 4.} {\it Under the conditions of existence of functions } $ ^r U^{\tau,\beta} (a,c;\alpha;\gamma;\delta;x) $ {\it such functional relations:}

\begin{gather}\label{eq42}
\frac{{\rm d}^n }{{\rm d} x^n}{} ^r U^{\tau,\beta} (a,c;\alpha;\gamma;\delta;x)=(-1)^n\frac{\Gamma(a+n)}{\Gamma(a)}{}^r U^{\tau,\beta} (a+n,c+n;\alpha;\gamma;\delta;x),
\end{gather}
\begin{gather}\label{eq43}
\frac{{\rm d}^n }{{\rm d} x^n} \left( e^{-x} \, {} ^r U^{\tau,\beta} (a,c;\alpha;\gamma;\delta;x) \right)=(-1)^n e^{-x} \, {}^r U^{\tau,\beta} (a,c+n;\alpha;\gamma;\delta;x),
\end{gather}
\begin{gather}\label{eq44}
^r U^{\tau,\beta} (a,c;\alpha+1;\gamma;\delta;x)= {}^r U^{\tau,\beta} (a,c;\alpha;\gamma;\delta;x)- \nonumber
\\
- \frac{r\tau}{\alpha}\frac{\Gamma(a-\delta)\Gamma(\alpha+\tau)\Gamma(\gamma)}{\Gamma(a)\Gamma(\alpha)\Gamma(\gamma+\beta)} {}^r U^{\tau,\beta} (a-\delta,c-\delta;\alpha+\tau;\gamma+\beta;\delta;x),
\end{gather}
\begin{gather}\label{eq45}
^r U^{\tau-\alpha,\beta-\gamma} (a,c;\alpha+1;\gamma;\delta;x)= {}^r U^{\tau-\alpha,\beta-\gamma} (a,c;\alpha;\gamma;\delta;x) - \nonumber
\\
- \frac{r(\tau - \alpha)}{\alpha}\frac{\Gamma(a-\delta)\Gamma(\tau)\Gamma(\gamma)}{\Gamma(a)\Gamma(\alpha)\Gamma(\beta)}{}^r U^{\tau-\alpha,\beta-\gamma} (a-\delta,c-\delta;\tau;\beta;\delta;x).
\end{gather}

\textbf{Proof.} We perform these functional relations (44), (45) by means of transformations $\alpha, \tau = \tau_1 + \alpha, \beta = \beta_1 + \gamma $, which are simple but cumbersome.

In this formula \eqref{eq39}, replace $\alpha +1$ for $\alpha$. Then we get:
\begin{gather}\label{eq46}
{} ^{r} _1 \Phi_1^{\tau,\beta} \left(\alpha +1;\gamma;- \frac{r} {t^{\delta}}\right)=\frac{\Gamma(\gamma)}{\Gamma(\alpha+1)}\sum_{n=0}^\infty\ \frac{\Gamma(\alpha+1+\tau n)}{\Gamma(\gamma+\beta n)}\left (-\frac{r}{t^{\delta}} \right) ^n \frac{1}{n!}.
\end{gather}
From formulas (46) and properties of Gamma-function \cite{Beitmen1965}: $\Gamma(\alpha +1)=\alpha \Gamma(\alpha)$, we get:
\begin{gather}
{} ^{r} _1 \Phi_1^{\tau,\beta} \left(\alpha +1;\gamma;- \frac{r} {t^{\delta}}\right)=\frac{\Gamma(\gamma)}{\alpha\Gamma(\alpha)}\sum_{n=0}^\infty\ \frac{(\alpha+\tau n)\Gamma(\alpha+\tau n)}{\Gamma(\gamma+\beta n)}\left (-\frac{r}{t^{\delta}} \right) ^n \frac{1}{n!}= \nonumber
\\
=\frac{\Gamma(\gamma)}{\Gamma(\alpha)}\sum_{n=0}^\infty\ \frac{\Gamma(\alpha+\tau n)}{\Gamma(\gamma+\beta n)}\left (-\frac{r}{t^{\delta}} \right) ^n \frac{1}{n!}+ \frac{\tau\Gamma(\gamma)}{\alpha\Gamma(\alpha)}\sum_{n=0}^\infty\ \frac{n\Gamma(\alpha+\tau n)}{\Gamma(\gamma+\beta n)}\left (-\frac{r}{t^{\delta}} \right) ^n \frac{1}{n!}= \nonumber
\\
={} ^{r} _1 \Phi_1^{\tau,\beta} \left(\alpha;\gamma;- \frac{r} {t^{\delta}}\right) + \frac{\tau\Gamma(\gamma)}{\alpha\Gamma(\alpha)}\sum_{n=1}^\infty\ \frac{n\Gamma(\alpha+\tau n)}{\Gamma(\gamma+\beta n)}\left (-\frac{r}{t^{\delta}} \right) ^n \frac{1}{n!}= \nonumber
\\
={} ^{r} _1 \Phi_1^{\tau,\beta} \left(\alpha;\gamma;- \frac{r} {t^{\delta}}\right) + \frac{\tau\Gamma(\gamma)}{\alpha\Gamma(\alpha)}\sum_{n=1}^\infty\ \frac{\Gamma(\alpha+\tau n)}{\Gamma(\gamma+\beta n)}\left (-\frac{r}{t^{\delta}} \right) ^n \frac{1}{(n-1)!}= \nonumber
\\
={} ^{r} _1 \Phi_1^{\tau,\beta} \left(\alpha;\gamma;- \frac{r} {t^{\delta}}\right) - \frac{r\tau\Gamma(\gamma)}{t^{\delta}\alpha\Gamma(\alpha)}\sum_{n=0}^\infty\ \frac{\Gamma((\alpha+\tau)+\tau n)}{\Gamma((\gamma+\beta)+\beta n)}\left (-\frac{r}{t^{\delta}} \right) ^n \frac{1}{n!}= \nonumber
\\
={} ^{r} _1 \Phi_1^{\tau,\beta} \left(\alpha;\gamma;- \frac{r} {t^{\delta}}\right) - \frac{r\tau\Gamma(\gamma)\Gamma(\alpha+\tau)}{t^{\delta}\alpha\Gamma(\alpha)\Gamma(\gamma+\beta)}{} ^{r} _1 \Phi_1^{\tau,\beta} \left(\alpha+\tau;\gamma+\beta;- \frac{r} {t^{\delta}}\right).
\end{gather}
Substituting the obtained value (47) into the generalization functions as Tricomi \eqref{eq38}, we get:
\begin{gather}
^r U^{\tau,\beta} (a,c;\alpha+1;\gamma;\delta;x)= {}^r U^{\tau,\beta} (a,c;\alpha;\gamma;\delta;x)- \nonumber
\\
- \frac{r\tau}{\alpha}\frac{\Gamma(\alpha+\tau)\Gamma(\gamma)}{\Gamma(a)\Gamma(\alpha)\Gamma(\gamma+\beta)} \int_{0}^\infty\ t^{a-\delta-1}(1+t)^{c-a-1}e^{-xt} \, {} ^{r} _1 \Phi_1^{\tau,\beta} \left(\alpha+\tau;\gamma+\beta;- \frac{r} {t^{\delta}}\right) {\rm d}t= \nonumber
\\
={}^r U^{\tau,\beta} (a,c;\alpha;\gamma;\delta;x)- \frac{r\tau}{\alpha}\frac{\Gamma(a-\delta)\Gamma(\alpha+\tau)\Gamma(\gamma)}{\Gamma(a)\Gamma(\alpha)\Gamma(\gamma+\beta)} {}^r U^{\tau,\beta} (a-\delta,c-\delta;\alpha+\tau;\gamma+\beta;\delta;x). \nonumber
\end{gather}
Thus, the property (44) is proved.

Through in (44) following replacement: $\tau=\tau_1 - \alpha$, $\beta=\beta_1 - \gamma$, we get formulas (45).

\end{document}